\documentclass{amsart}
\usepackage{amssymb,amsmath}
\usepackage{ytableau}
\usepackage[enableskew]{youngtab}
\usepackage{young}
\usepackage{tikz,pgf}

\newtheorem{theorem}{Theorem}

\newtheorem{lemma}{Lemma}

\begin{document}

\title{Euler's partition theorem for all moduli and new companions to Rogers-Ramanujan-Andrews-Gordon identities}
\author{Xinhua Xiong, William J. Keith}
\maketitle

\begin{abstract}
We generalise Euler's partition theorem involving odd parts and different parts for all moduli and provide new companions to Rogers-Ramanujan-Andrews-Gordon identities related to this theorem.
\end{abstract}

\section{Introduction}
In the theory of partitions, Euler's partition theorem involving odd parts and different parts is one of the famous theorems. It claims that the number of partitions of $n$ into odd parts is equal to  the number of partitions $n$ into different parts. By constructing a bijection, Sylvester \cite{Syl82} not only proved Euler's theorem, but also provided a refinement of it which can be stated as, ``the number of partitions of $n$ into odd parts with exactly $k$ different parts is equal to the number of partitions of $n$ into different parts such that exactly $k$ sequences of consecutive integers occur in each partition.''  Bessenrodt  \cite{bes94} proved that Sylvester's bijection implies that the number of partitions of $n$ into different parts with the alternating sum $\Sigma$ is equal to the number of partitions of $n$ with $\Sigma$ odd parts. Here for a partition $\lambda=(\lambda_1,\lambda_2,\dots,\lambda_{k})$, the alternating sum is defined by 
\begin{equation}
\Sigma=\lambda_{1}-\lambda_{2}+\lambda_{3}-\lambda_{4}+\dots+(-1)^{k+1}\lambda_{k}.\label{alternating}
\end{equation}
 In \cite{Kim-Yee}, Kim and Yee gave a different description of Sylvester's bijection which provides a simpler proof of the refinement of Euler's theorem due to Bessenrodt. There are several other refinements and variants of Euler's theorem. See \cite{All3,AndrewsEuler,Ber1,Ber2, Ber3, Pak, Yee2001,Yee2002}.

We can think of Euler's theorem as a theorem on partitions involving modulus two by interpreting odd parts as parts $\equiv 1\pmod 2.$  The first nontrivial generalisation of Euler's theorem for all moduli in this sense is the following theorem due to Pak-Postnikov.

 \begin{theorem}[Pak-Postnikov \cite{Pak-Postnikov}]   \label{puretype thm} 
 The number of partitions of $n$ with type $(c, m-c, c, m-c,\dots)$  is equal to the number of partitions of $n$ with all parts $\equiv c \pmod m$.
 \end{theorem}
By the type $(c, m-c,c,m-c,\dots)$ for a partition $\lambda$, it means that $\lambda$ has the length divisible by $m$ by allowing zero as parts  and has $c\geq 1$ largest parts, $m-c$ second largest parts, etc. So it has the form:
\begin{gather*}
\lambda_{1}=\lambda_{2}=\lambda_{3}=\dots=\lambda_{c}>\lambda_{c+1}=\lambda_{c+2}=\\
 \dots = \lambda_{m}>\lambda_{m+1}=\dots=\lambda_{m+c}>\dots
\end{gather*}
The second author's doctoral thesis \cite{KeithThesis} has a chapter devoted to various identities of this nature, such as for $m$-\emph{falling} or $m$-\emph{rising} partitions, those in which least positive residues of each part modulo $m$ form a nonincreasing or nondecreasing sequence.

In this paper, we prove a theorem which generalises Euler's partition theorem mentioned above for all moduli, and simultaneously generalises Pak-Postnikov's theorem. We also show that this theorem provides new companions to Rogers-Ramanujan-Andrews-Gordon identities. 
 
In Section 2 we give all definitions necessary to state our theorem.  In Section 3 we prove the theorem by a bijection originally due to Stockhofe in \cite{Stockhofe} (an English translation of the original German can be found as an appendix to \cite{KeithThesis}), slightly extended and heavily specialized for the present purpose.  The original map was a general bijection on all partitions; we will show that the properties required hold when specialized to the sets of interest for the theorem.  All concepts and claims necessary will be defined and proved here to keep the paper self-contained.

\section{Statement of theorem}
In order to state our theorem, we introduce some notations and terminologies.

$\bullet$  For a partition $\lambda=(\lambda_1,\lambda_2, \lambda_{3},\dots,\lambda_{r})$ of $n$, we denote $n$ by $|\lambda |.$ We sometimes write a partition as the form
$$\lambda_1+\lambda_2+\lambda_3+\dots+\lambda_{r}  \quad \mbox{\quad or the form\quad} \quad   \lambda_1\geq\lambda_2\geq\lambda_{3}\geq\dots\geq\lambda_{r}.$$

$\bullet$ Let $m\geq 2$ be an integer, which is the moduli in our sense.  For a partition $\lambda_1\geq\lambda_2\geq\lambda_{3}\geq\dots\geq \lambda_{km}$ ($k\geq1$ ) with length divisible by $m$, we define its alternating sum type to be an $(m-1)$-tuple of non-negative integers $(\Sigma_{1},  \Sigma_{2},  \dots, \Sigma_{m-2},  \Sigma_{m-1})$ given by
\begin{gather*}
\Sigma_{1}=\sum_{i=1}^{k}\lambda_{(i-1)m+1}-\lambda_{(i-1)m+2},\\
 \Sigma_{2}=\sum_{i=1}^{k}\lambda_{(i-1)m+2}-\lambda_{(i-1)m+3},\\
   \cdots\\
    \Sigma_{m-1}=\sum_{i=1}^{k}\lambda_{(i-1)m+m-1}-\lambda_{im}.
\end{gather*}

For example, if $m=3$, then the partition $6+5+4+3+2+1$ has the alternating sum type $(6-5+3-2,5-4+2-1)=(2, 2)$. The partitions appearing in the theorem of Pak-Postnikov have the special alternating sum type $(0,0,\dots,0, \Sigma, 0,\dots,0)$, where $\Sigma$ lies at the $c^{th}$ position. If $m=2$, the alternating sum type is exactly the alternating sum is given by (\ref{alternating}).


$\bullet$ When we speak of alternating sum types for a partition of length $k$, we allow the last few parts to be zero so that any partition has length $\lceil k/m \rceil m$. For example if $m=3$, the partition $5+4+3+3$ has length $6$ by viewing it as $5+4+3+3+0+0$ and it has two basic units: $5+4+3$ and $3+0+0.$

$\bullet$ For a partition $\lambda$ with all parts $\not\equiv 0\pmod m$, known as an $m$-regular partition, we define its length type to be the $(m-1)$-tuple of non-negative $(l_{1}, l_{2}, l_{3}, \dots, l_{m-2},  l_{m-1})$, where for $1\leq i\leq m-1$,  $l_{i}$ is the number of parts of $\lambda$ which are congruent to $i$ modulo $m$. For example, for a partition $\lambda$ with all parts congruent to $c$, its length type is $(0,0,\dots,0,l,0\dots,0)$, where $l$ is the number of parts of $\lambda$ and lies at the $c^{th}$ position.

$\bullet$ Given a partition $\lambda = (\lambda_1, \dots, \lambda_k)$, its conjugate $\lambda^\prime$ is $(\#\{\lambda_i \geq 1\},\#\{\lambda_i \geq 2\},\dots)$.


Now we can state our theorem.
\begin{theorem}
Let $m\geq 2$. Let $P$ be the set of partitions in which each part can be repeated at most $m-1$ times. (This implies that their alternating sum types cannot be $(0,0,\dots,0)$.) Let $Q$ be the set of partitions with no parts $\equiv 0\pmod m$. Then we have the partition identity:
$$
\sum_{\lambda \in P} z_{1}^{\Sigma_{1}(\lambda)}  z_{2}^{\Sigma_{2}(\lambda)      }  \dots z_{m-1} ^{\Sigma_{m-1}(\lambda)          }  q^{|\lambda |}      =        \sum_{\mu \in Q} z_{1}^{l_{1}(\mu)}  z_{2}^{l_{2}(\mu)      }  \dots z_{m-1} ^{l_{m-1}(\mu)          }  q^{|\mu |}.   
 $$
Equivalently,  the number of partitions  of $n$ with the alternating sum type $(\Sigma_{1}, \Sigma_{2}, \dots \Sigma_{m-1})$ is equal to the number of partitions of $n$ with $\Sigma_{1}$ parts congruent to $1$ modulo $m$, $\Sigma_{2}$ parts congruent to $2$ modulo $m$, $\dots$ , $\Sigma_{m-1}$ parts congruent $m-1$ modulo $m$.
\end{theorem}

If we let $z_{1}=z_{2}=\dots z_{m-1}=z$, we get the result that the number of partitions of $n$ with parts repeated at most $m-1$ times and total alternating sum $\Sigma_{1}+\Sigma_{2}+\dots +\Sigma_{m-1}$ is equal to the number of partitions of $n$ with no parts congruent to $0$ modulo $m$  and $\Sigma_{1}+\Sigma_{2}+\dots +\Sigma_{m-1}$ parts, which is a refinement of Glaisher's theorem:
\begin{theorem}[Glaisher \cite{Pak}]
The number of partitions of $n$ with parts repeated  at most $m-1$ times is equal to the number of partitions of $n$ with no parts is congruent to $0$ modulo $m$.
\end{theorem}
When the alternating sum type is pure type, this theorem reduces to Theorem 1.1 due to Pak-Postnikov. When $m$ is $2$, this theorem reduces to the refinement of Euler's theorem due to Bessenrodt, Kim and Yee.

We give $n=10$ and $m=3, 4$ to illustrate this theorem. 
We list partitions in $P$, their alternating sum types and the numbers on the left, and the corresponding parts for partitions in $Q$ on the right. We only list all partitions with mixed types.

\begin{minipage}{0.38\textwidth}
\small
\begin{align*}
& \mbox{ Partitions in $P$} &(\Sigma_{1}, \Sigma_{2})&\quad\sharp \\
&\left\{
\begin{array}{llll}
&3+3+2+2+1\\
&5+4+2\\
&4+4+2+1\\
&4+3+2+1+1
\end{array}
\right\} &(1, 2)   &\quad4 \\
&\left\{
\begin{array}{llll}
&6+3+2\\
&5+3+2+1\\
&5+2+2+1+1\\
&4+3+2+2
\end{array}
\right\} &(3,1)  &\quad4\\
&\left\{
\begin{array}{ll}
&6+4+1\\
&5+4+1+1
\end{array}
\right\} &(2,3)  &\quad2 \\
&\left\{
\begin{array}{ll}
&7+3+1\\
&6+3+1+1
\end{array}
\right\} &(4,2)  &\quad2\\
&\left\{
\begin{array}{ll}
&8+2+1\\
&7+2+1+1
\end{array}
\right\} &(6,1)  &\quad2\\
&\left\{
\begin{array}{l}
6+5
\end{array}
\right\} &(1,5)  &\quad1\\
&\left\{
\begin{array}{l}
7+4
\end{array}
\right\} &(3,4)  &\quad1\\
&\left\{
\begin{array}{l}
8+3
\end{array}
\right\} &(5,3)  &\quad1\\
&\left\{
\begin{array}{l}
9+2
\end{array}
\right\} &(7,2)  &\quad1\\
&\left\{
\begin{array}{l}
10+1
\end{array}
\right\} &(9,1)  &\quad1
 \end{align*}
\end{minipage}
\begin{minipage}{0.50\textwidth}
\small
\begin{align*}
& \mbox{ Partitions in $Q$ } &(l_{1}, l_{2})  &\quad\sharp\\
&\left\{
\begin{array}{llll}
&8+2+1\\
&7+2+2\\
&5+5+1\\
&5+4+2
\end{array}
\right\} &(1,2)  &\quad4\\
&\left\{
\begin{array}{llll}
&8+1+1+1\\
&7+2+1+1\\
&5+4+1+1\\
&4+4+2+1
\end{array}
\right\} &(3,1)  &\quad4\\
&\left\{
\begin{array}{ll}
&5+2+2+1+1\\
&4+2+2+2+1
\end{array}
\right\} &(2,3)  &\quad2\\
&\left\{
\begin{array}{ll}
&5+2+1+1+1+1\\
&4+2+2+1+1+1
\end{array}
\right\} &(4,2)  &\quad2\\
&\left\{
\begin{array}{ll}
&5+1+1+1+1+1+1\\
&4+2+1+1+1+1+1
\end{array}
\right\} &(6,1)  &\quad2\\
&\left\{
\begin{array}{l}
2+2+2+2+2+1
\end{array}
\right\} &(1,5)  &\quad1\\
&\left\{
\begin{array}{l}
2+2+2+2+1+1+1
\end{array}
\right\} &(3,4)  &\quad1\\
&\left\{
\begin{array}{l}
2+2+2+1+1+1+1+1
\end{array}
\right\} &(5,3)  &\quad1\\
&\left\{
\begin{array}{l}
2+2+1+1+1+1+1+1+1
\end{array}
\right\} &(7,2)  &\quad1\\
&\left\{
\begin{array}{l}
2+1+1+1+1+1+1+1+1+1
\end{array}
\right\}&(9,1)  &\quad1
\end{align*}  
\end{minipage}

\vspace{0.5cm}



\begin{minipage}{0.43\textwidth}
\small
\begin{align*}
\small
& \mbox{ Partitions in $P$} &(\Sigma_{1}, \Sigma_{2},\Sigma_{3})  &\sharp \\
&\left\{
\begin{array}{lll}
&4+3+2+1\\
&3+3+2+1+1\\
&3+2+2+1+1+1
\end{array}
\right\}\, &(1,1,1) \,&3\\
&\left\{
\begin{array}{ll}
&5+3+1+1\\
&4+3+1+1+1
\end{array}
\right\}\, &(2,2,0) \, &2\\
&\left\{
\begin{array}{ll}
&5+2+2+1\\
&4+2+2+1+1
\end{array}
\right\}\, &(3,0,1) \, &2\\
&\left\{
\begin{array}{ll}
&6+2+1+1\\
&5+2+1+1+1
\end{array}
\right\}\, &(4,1,0) \, &2
\end{align*}
\end{minipage}
\begin{minipage}{0.50\textwidth}
\small
\begin{align*}
& \mbox{ Partitions in $Q$ } &(l_{1}, l_{2},l_{3}) \, &\sharp\\
&\left\{
\begin{array}{lll}
&5+3+2\\
&6+3+1\\
&7+2+1
\end{array}
\right\}\, &(1,1,1) \, &3\\
&\left\{
\begin{array}{ll}
&5+2+2+1\\
&6+2+1+1
\end{array}
\right\}\, &(2,2,0) \, &2\\
&\left\{
\begin{array}{ll}
&5+3+1+1\\
&7+1+1+1
\end{array}
\right\}\, &(3,0,1) \, &2\\
&\left\{
\begin{array}{ll}
&5+2+1+1+1\\
&6+1+1+1+1
\end{array}
\right\}\, &(4,1,0) \, &2
\end{align*}
\end{minipage}

\begin{minipage}{0.43\textwidth}
\small
\begin{align*}
&\left\{
\begin{array}{l}
4+4+2
\end{array}
\right\}\,\quad\quad\quad\quad\quad\quad\quad &(0,2,2)  &1\\
&\left\{
\begin{array}{l}
4+3+3
\end{array}
\right\}\, &(1,0,3) \, &1\\
&\left\{
\begin{array}{l}
5+4+1
\end{array}
\right\}\, &(1,3,1) \, &1\\
&\left\{
\begin{array}{l}
5+3+2
\end{array}
\right\}\, &(2,1,2) \, &1\\
&\left\{
\begin{array}{l}
6+4
\end{array}
\right\}\, &(2,4,0) \, &1\\
&\left\{
\begin{array}{l}
6+3+1
\end{array}
\right\}\, &(3,2,1) \, &1\\
&\left\{
\begin{array}{l}
6+2+2
\end{array}
\right\}\, &(4,0,2) \, &1\\
&\left\{
\begin{array}{l}
7+3
\end{array}
\right\}\, &(4,3,0) \, &1\\
&\left\{
\begin{array}{l}
7+2+1
\end{array}
\right\}\, &(5,1,1) \, &1\\
&\left\{
\begin{array}{l}
8+2
\end{array}
\right\}\, &(6,2,0) \,&1\\
&\left\{
\begin{array}{l}
8+1+1
\end{array}
\right\}\, &(7,0,1) \, &1\\
&\left\{
\begin{array}{l}
9+1
\end{array}
\right\}\, &(8,1,0) \, &1
\end{align*}
\end{minipage}
\begin{minipage}{0.50\textwidth}
\small
\begin{align*}
\vspace{0.5cm}
&\left\{
\begin{array}{l}
3+3+2+2
\end{array}
\right\}\, &(0,2,2) \, &1\\
&\left\{
\begin{array}{l}
3+3+3+1
\end{array}
\right\}\, &(1,0,3) \, &1\\
&\left\{
\begin{array}{l}
3+2+2+2+1
\end{array}
\right\}\, &(1,3,1) \, &1\\
&\left\{
\begin{array}{l}
3+3+2+1+1
\end{array}
\right\}\, &(2,1,2) \, &1\\
&\left\{
\begin{array}{l}
2+2+2+2+1+1
\end{array}
\right\}\, &(2,4,0) \, &1\\
&\left\{
\begin{array}{l}
3+2+2+1+1+1
\end{array}
\right\}\, &(3,2,1) \, &1\\
&\left\{
\begin{array}{l}
3+3+1+1+1+1
\end{array}
\right\}\, &(4,0,2) \, &1\\
&\left\{
\begin{array}{l}
2+2+2+1+1+1+1
\end{array}
\right\}\, &(4,3,0) \, &1\\
&\left\{
\begin{array}{l}
3+2+1+1+1+1+1
\end{array}
\right\}\, &(5,1,1) \, &1\\
&\left\{
\begin{array}{l}
2+2+1+1+1+1+1+1
\end{array}
\right\}\, &(6,2,0) \, &1\\
&\left\{
\begin{array}{l}
3+1+1+1+1+1+1+1
\end{array}
\right\}\, &(7,0,1) \, &1\\
&\left\{
\begin{array}{l}
2+1+1+1+1+1+1+1+1
\end{array}
\right\}\, &(8,1,0) \,&1
\end{align*}  
\end{minipage}



\section{Proof of the main theorem}

We begin with a simple lemma.

\begin{lemma}\label{conj} The conjugates $\lambda^{\prime}$ of partitions $\lambda$ with alternating sum type $$(s_1,\dots,s_{m-1})$$ are precisely those partitions of length type $(s_1,\dots,s_{m-1})$.
\end{lemma}

\emph{Proof.} Suppose $\lambda_{km+i} - \lambda_{km+i+1} = c$ contributes a nonzero amount to $s_i$.  Then in the conjugate partition, $c$ parts of size $km+i$ appear.  The converse also holds.

Call $m$-\emph{flat} a partition in which all differences between consecutive parts are strictly less than $m$ and the smallest part is less than $m$.  These are clearly the conjugates of partitions in $P$.  We will prove by bijection that

\begin{theorem} There is a bijection between $m$-regular partitions of any given length type $(\ell_1,\dots,\ell_{m-1})$ and $m$-flat partitions of the same length type.
\end{theorem}

\emph{Proof.} In fact, the bijection even preserves the sequential order of the nonzero residues modulo $m$; our map will consist of rearranging units of size $m$.

It is useful to define two operations analogous to scalar multiplication and vector addition on partitions.  For convenience, assume that all partitions are equipped with infinite tails consisting solely of zeros.  The scalar multiple of a partition $\lambda = (\lambda_1,\lambda_2,\dots)$ by the positive integer $m$ is the partition $m \lambda = (m \lambda_1,m \lambda_2,\dots)$.  Given two partitions $\lambda = (\lambda_1,\lambda_2,\dots)$ and $\mu = (\mu_1,\mu_2,\dots)$, we can define their (infinite-dimensional) vector sum $\lambda + \mu = (\lambda_1+\mu_1,\lambda_2 + \mu_2,\dots)$.  (Partition addition in the literature sometimes means taking the nonincreasing sequence of the multiset union of all parts of both partitions; we will not require this operation.)

Begin with an $m$-flat partition.  We will remove multiples of $m$ to construct a partition $\pi$, via a sequence of intermediate partitions $\lambda^{(0)}$, $\lambda^{(1)}$, $\lambda^{(2)}$, etc.  As we do so, we will use the removed multiples of $m$ to construct a second partition $m\sigma = (m\sigma_1,m\sigma_2,\dots)$.

Initialize $\sigma = ()$, the empty partition.  

\textbf{Step 1.} First, construct $\lambda^{(0)}$ by removing from $\lambda$ any parts divisible by $m$ for which, after removal, the partition $\lambda^{(0)}$ is still $m$-flat.  These will be parts such that

\begin{itemize}
\item $\lambda_i = k_i m = \lambda_{i+1}$, i.e. all but the last of a repeated part divisible by $m$; 
\item $\lambda_1 = k_1 m$ if the largest part is divisible by $m$; or
\item parts $\lambda_i = k_i m$, $i>1$, such that $\lambda_{i-1} = k_i m+j_1$, $\lambda_{i+1} = (k_i-1)m + j_2$, with $0 < j_1 < j_2 < m$.
\end{itemize}

Thus the remaining parts divisible by $m$ in $\lambda^{(0)}$ are all distinct, not the largest (or smallest) parts, and any remaining part $\lambda_i = k_i m$ lies between $\lambda_{i-1} = k_i m + j_1$ and $\lambda_{i+1} = (k_i-1)m+j_2$ with $0 < j_2 \leq j_1 < m$.

For each part $\lambda_i = k_i m$ removed, add $k_i m$ to $m\sigma$ as a part.

For the previous step, the order of removal did not matter, although of course $m \sigma$ is arranged in nonincreasing order.  In the next step, we work from the largest part divisible by $m$ to the smallest.

\textbf{Step 2.} Begin with $\lambda^{(0)}$ and set $j=0$.

\begin{enumerate}
\item If $\lambda^{(j)}$ has no parts divisible by $m$, stop.
\item If $\lambda_i = k_i m$ is the largest part in $\lambda^{(j)}$ divisible by $m$, remove $\lambda_i$ from $\lambda^{(j)}$.  Renumber following parts.
\item In addition, reduce by $m$ all parts $\lambda_1$ through $\lambda_{i-1}$.  The remaining partition is now $\lambda^{(j+1)}$.  Increment $j$.
\item Add $m(k_i + i - 1)$ to $m\sigma$ as a part.
\item Repeat.
\end{enumerate} 

The following lemma concerning parts removed in Step 2 will be useful when we wish to prove that this process is reversible.

\begin{lemma}\label{AngleLemma} Parts added to $\sigma$ in Step 2 are always at least the size of those removed in Step 1, and are added in nondecreasing order of size.  The largest possible size of a part added to $m\sigma$ in Step 2 is the number of parts in $\lambda$ not divisible by $m$.
\end{lemma}

\emph{Proof.} If $\lambda$ is an $m$-flat partition and $\lambda_i = k_im+j_1$, $j_1 \not\equiv 0 \pmod{m}$, with $\lambda_{i+c} = (k_i-1)m+j_2$ the next smaller part which is nonzero modulo $m$, this necessarily requires $0 < j_1 < j_2 < m$.  In this case refer to $\lambda_i$ as a \emph{descent} of $\lambda$.  If a part $k_i m$ appears between $\lambda_i$ and $\lambda_{i+c}$ defining a descent, we speak of $k_i m$ as appearing within the descent.

Parts divisible by $m$ are not the largest part of $\lambda^{(j)}$ and do not appear within descents of $\lambda^{(j)}$: these were removed in Step 1.

Suppose part $\lambda_i = k_im$ appears in $\lambda^{(j)}$, so that when removed we will add part $m(k_i+i-1)$ to $m\sigma$.  The next part, if any, which will be removed is $\lambda_{i+s-1} = k_{i+s-1}m$, $s\geq 2$, after the renumbering.  (That is, it was $\lambda_{i+s}$ before renumbering.)

We have $k_{i+s-1} < k_i$, decreasing by exactly 1 for each descent passed as we read from (after renumbering) $\lambda_i$ to $\lambda_{i+s-1}$, plus 1 immediately, in essence thinking of the passage from $k_im$ to $\lambda_i$ as a descent.  The total decrease is at most $s-1$, since $\lambda_{i+s-2}$ and $\lambda_{i+s-1}$ cannot be descents (parts $k_i m$ in Step 2 do not appear within descents).

On the other hand, the number of parts added due to subtraction from previous parts always increases by $s-1$: one for each part passed regardless of whether it is a descent or not, less 1 because 1 fewer part exists prior to $\lambda_{i+s-1}$ after removal of $\lambda_i$.  Thus we have $k_{i+s-1}+(i+s-1)-1 \geq k_i+i-1$.

This likewise holds for the first part removed in Step 2, taking $i=0$ for a potential largest part.  The first $k_i+i-1$ removed in Step 2 is decreased from this size by at most $i$ and increased by $i$ exactly.

Finally, the largest a part removed in Step 2 can be is if we add as much as possible, with steps across descents being irrelevant; that is, the largest possible part that could be removed is a part $\lambda_i = m$ which is the next-to-last part, followed by a single part not divisible by $m$.  Since all previous parts divisible by $m$ would have been removed at this step, clearly in this case $1+i-1$ is exactly the number of parts in $\lambda$ not divisible by $m$.  By the previous clauses, this is the largest removal.

Thus all claims of the lemma hold.

\textbf{Step 3.} After Step 2, we now have some $\lambda^{(j)}$ which is simultaneously $m$-regular and $m$-flat, and $m\sigma$ consisting of parts divisible by $m$. Set $\pi = \lambda^{(j)}$.  Our final partition is $\pi+m(\sigma^{\prime}).$  Since by our lemma the largest part of $\sigma$ was less than or equal to the number of parts in $\lambda$ not divisible by $m$, its conjugate has at most this number of parts, so we only add multiples of $m$ to such parts in $\pi$.  The resulting partition has all parts not divisible by $m$.

Since no step in this construction alters the residue modulo $m$ of a part not divisible by $m$, it is an easy lemma that

\begin{lemma} The length type of the parts of $\lambda$ not divisible by $m$, read as a partition, is the same as the length type of the partition $\pi + m (\sigma^\prime)$.
\end{lemma}

We now briefly show that the map is reversible.  Starting with a partition $\mu$ into parts not divisible by $m$, we will construct a sequence of partitions $\pi^{(j)}$ which begin with the $m$-flat, $m$-regular portion of $\mu$ and have parts from $\sigma$ inserted.

\textbf{Step 3 Reverse.} It is easy to break a partition $\mu$ into a flat plart plus a collection of parts divisible by $m$, as $\pi + m(\sigma^{\prime})$.  Whenever $\mu_i - \mu_{i+1} \geq m$ (including for the smallest part: treat the next part as 0), add 1 to $\sigma^{\prime}$ for parts 1 through $i$.  Subtract $m$ from all parts $\mu_1$ through $\mu_i$.  Repeat.  When done with all possible removals, conjugate $\sigma^{\prime}$ to obtain $\sigma$.

\textbf{Step 2 Reverse.}  Observe that if we wish to insert a part $m\sigma_1$ into $\pi$, we must determine whether it is to be inserted in the reverse of Step 2 or Step 1.  Step 2 insertions occur when part $m\sigma_1$ is larger than $k_1$ for $(\pi^{(j)})_1 = k_1 m + j_1$.  The position where such a part can be inserted is unique, since by the proof of Lemma \ref{AngleLemma} there can be only one $i$ such that $\sigma_1 = k_i + i - 1$ and in which $m \sigma_1$ would not be appearing within a descent.  Passing a column that is not a descent changes the amount to be added; passing a column that is a descent does not, but is not a place where parts are added in Step 2.

\textbf{Step 1 Reverse.} This step is easy since the order in which parts are inserted will not matter.  Once parts are small enough that they can be inserted into $\pi^{(j)}$ while retaining flatness, insert all at once.  A part of size $k_i m$ will go precisely after a part of size $k_i m$ if one already exists, or within the descent at $\lambda_i = k_i m+j_1$ if the next part is not divisible by $m$.

The result is our desired $m$-flat partition.

An example of the bijection may be illustrative.  Let our modulus be $m=5$.

Let our starting partition be $(9,9,8,8,8,7,6,6,6,6,5,5,5,4,4,2,2,2,2,1,1,1)$.  We observe that its alternating sum type is $(1,2,1,1)$.  Its conjugate is $\lambda = (22,19, 15,15, 13,10,6,5,2)$.  

Write the 5-modular diagram of this partition:

\begin{center}
\begin{tabular}{ccccc}
2 & 5 & 5 & 5 & 5 \\
4 & 5 & 5 & 5 & \\
5 & 5 & 5 & & \\
5 & 5 & 5 & & \\
3 & 5 & 5 & & \\
5 & 5 & & & \\
1 & 5 & & & \\
5 & & & & \\
2 & & & & 
\end{tabular}
\end{center}

The first parts we remove are those which can be removed whole without destroying 5-flatness.  We can remove the parts 5 (because $6-2=4$), not the 10, and one of the 15s but not the second (the first because it is repeated, but not the second because $19-13 = 6$).  So far $\sigma \cdot 5 = (15, 5) = (3,1)\cdot 5$.

The remaining partition is now $(22,19, 15,13,10,6,2)$.

\begin{center}
\begin{tabular}{ccccc}
2 & 5 & 5 & 5 & 5 \\
4 & 5 & 5 & 5 & \\
5 & 5 & 5 & & \\
3 & 5 & 5 & & \\
5 & 5 & & & \\
1 & 5 & & & \\
2 & & & & 
\end{tabular}
\end{center}

We remove the part 15, and in addition subtract 5 from 22 and 19.  Thus, we add the part 25 to $5\sigma$, obtaining $(25,15,5)$ so far, and are left with the following partition, $(17, 14, 13, 10, 6, 2)$:

\begin{center}
\begin{tabular}{ccccc}
2 & 5 & 5 & 5 & \\
4 & 5 & 5 & & \\
3 & 5 & 5 & & \\
5 & 5 & & & \\
1 & 5 & & & \\
2 & & & & 
\end{tabular}
\end{center}

Finally we remove the 10 and a 5 from each of the previous three larger parts, adding a 25 to $\sigma$ and finishing with $\sigma \cdot 5 = (25,25,15,5) = (5,5,3,1)\cdot 5$, and $\pi = (12,9,8,6,2)$:

\begin{center}
\begin{tabular}{ccccc}
2 & 5 & 5 & & \\
4 & 5 & & & \\
3 & 5 & & & \\
1 & 5 & & & \\
2 & & & & 
\end{tabular}
\end{center}

To combine these into a new partition we conjugate $\sigma$, obtaining $\sigma^{\prime} = (4,3,3,2,2)$, and add 5 times this partwise to $\pi$:

\begin{center}
\begin{tabular}{ccccccc}
2 & 5 & 5 & 5 & 5 & 5 & 5 \\
4 & 5 & 5 & 5 & 5 & & \\
3 & 5 & 5 & 5 & 5 & & \\
1 & 5 & 5 & 5 & & & \\
2 & 5 & 5 & & & & 
\end{tabular}
\end{center}

Our final partition is $(32, 24, 23, 16, 12)$.  Its length type is $(1,2,1,1)$, as desired.

If we were to reverse our map, we would observe that $\sigma^{\prime}$ has five parts of size at least 2 since the smallest part of $\mu$ is $2+5+5$, and so forth obtain $\sigma$; observing that the largest part of $\sigma$ is 5, we determine that we should set $i=3$, since setting $i=4$ is too large (adding part 5 and following it by $1+5+5$ would not result in a flat partition), whereas $i=2$ would not result in a partition at all (15 preceded by 14).  The other insertions are likewise unique.

\section{New companions to Rogers-Ramanujan-Andrews-Gordon identities}

Besides Euler's partition theorem involving odd parts and different parts, Rogers-Ramanujan-Andrews-Gordon identities are another famous partition theorem; see \cite{AlladiBerkovich2, AndrewsRR, Ber-peter, Bressoud}. Recall that the first Rogers-Ramanujan identity (partition version) says that the number of partitions of $n$ with the condition that the difference between any two parts is at least $2$ (called Rogers-Ramanujan partitions) is equal to the number of partitions of $n$ such that each part is congruent to $1$ or $4$ modulo $5$. From our viewpoint, partitions with each part congruent to $1$ or $4$ modulo $5$ are exactly partitions belonging to $Q$ with length types $(l_{1},0, 0, l_{4})$, $(l_{1},l_{4})\neq (0,0)$. Then our theorem gives the following companion to the first Rogers-Ramanujan identity:
\begin{theorem}
The number of partitions of $n$ where the difference between any two parts is at least
$2$ is equal to the number of partitions of $n$ with parts repeated  at most $4$ times and alternating sum types $(\Sigma_{1},0,0,\Sigma_{4})$,  where $(\Sigma_{1}, \Sigma_{4})\neq (0,0)$.
\end{theorem}
We give an example to illustrate this theorem. Let $n=11$, the partitions of $11$ with the condition that the difference is at least
$2$ are
\begin{align*}
11,\quad 10+1,\quad 9+2,\quad 8+3, \quad
7+4, \quad 7+3+1,\quad 6+4+1.
\end{align*}
And the partitions of $11$ with alternating sum types $(\Sigma_{1},0,0,\Sigma_{4})$ are 
$$3+2+2+2+2\,(1,0,0,0),\quad \,3+2+2+2+1+1\,(2,0,0,1),$$
$$\quad 4+2+2+2+1\,(2,0,0,1), \quad\,7+1+1+1+1\,(6,0,0,0),$$ 
$$\quad    5+2+2+2\,(3,0,0,2),\quad 8+1+1+1\,(7,0,0,1),  \quad 11\,(11,0,0,0).$$
We list the alternating sum type following each partition.
We have a similar companion on the second Rogers-Ramanujan identity:
\begin{theorem}
The number of partitions of $n$ where the difference between any two parts is at least
$2$ and $1$ is not a part is equal to the number of partitions of $n$ with parts repeated at most $4$ times and alternating sum type $(0,\Sigma_{2},\Sigma_{3},0)$ and $(\Sigma_{2}, \Sigma_{3})\neq (0,0)$.
\end{theorem}
We still use $n=11$ to illustrate this theorem.  The partitions of $11$ with the condition that the difference is at least
$2$ and $1$ is not a part are
\begin{align*}
11,\quad 9+2,\quad 8+3,\quad
7+4.
\end{align*}
And the partitions of $11$ with alternating sum types $(0,\Sigma_{2},\Sigma_{3},0)$ are 
$$3+3+3+1+1\,(0,0,2,0),\quad\quad 4+4+1+1+1\,(0,3,0,0),$$
$$ 4+4+3\,(0,1,3,0)\quad\quad 5+5+1\,(0,4,1,0).$$
We list the corresponding alternating sum type following each partition.

For Andrews-Gordon's identities, we have 
\begin{theorem}
Let $d\geq 1$, $1\leq i\leq 2d$.
Suppose the conjecture is true, then the number of partitions $\lambda_{1}+\lambda_{2}+\lambda_{3}+\dots+\lambda_{r}$ of $n$ such that no more than $i-1$ of the parts are $1$ and pairs of consecutive integers appear at most $d-1$
 times is equal to the number of partitions of $n$ with parts repeated  at most $2d$ times and alternating sum type $(\Sigma_{1},\Sigma_{2},\dots,\Sigma_{2d-1},\Sigma_{2d})\neq(0,0,\dots,0)$ satisfying that both $\Sigma_{i}$ and $\Sigma_{2d+1-i}$ are zero.
\end{theorem}

\section*{Acknowledgements}
The first author would like to thank Professor Peter Paule and Professor Christian Krattenthaler for their  comments on an earlier version of this paper during the Strobl meeting. The first author was supported  by the Austria Science Foundation (FWF) grant SFB F50-06 (Special Research Program ``Algorithmic and Enumerative Combinatorics'') and partially supported by the National Natural Science Foundation of China (11101238).

\end{document}